\documentclass[11pt,reqno,a4paper,twoside]{amsart}
\usepackage{amssymb,amsmath,amscd,amstext,amsthm,amsfonts}
\usepackage[ansinew]{inputenc} 
\usepackage{graphicx}
\usepackage[mathscr]{eucal}

\numberwithin{equation}{section}

\newcommand{\R}{\mathbb{R}}
\newcommand{\Pp}{\mathbb{P}}
\newcommand{\C}{\mathbb{C}}

\newcommand{\Dsk}{\mathbb{D}}

\newcommand{\N}{\mathbb{N}}
\newcommand{\Z}{\mathbb{Z}}

\newcommand{\Su}{\mathbb{S}}

\newcommand{\obar}[1]{{\overline{#1}}}
\newcommand{\OO}{{\rm O}}

\newcommand{\SL}{{\rm SL}}
\newcommand{\PSL}{{\rm PSL}}

\newcommand{\matriz}[4]{ 
\left[
\begin{array}{cc}
{#1}&{#2}\\
{#3}&{#4}
\end{array}\right] }
\newcommand{\Mod}[1]{\left\vert{#1}\right\vert}
\newcommand{\nrm}[1]{\left\|#1\right\|}

\newcommand{\figura}[2]{
  \begin{center}
 \includegraphics*[scale={#2}]{{#1}} 
 \end{center}
}

\newtheorem{lema}{Lemma}

\newtheorem{prop}{Proposition}
\newtheorem{teor}{Theorem}

\newcommand{\dem}{ \par\medbreak\noindent{\bf
Proof. }\enspace} 
 
\newcommand{\cqd}{\;
$\sqcup\!\!\!\!\sqcap\bigskip$}

\begin{document}

\title[Herman-Avila-Bochi formula]{A new proof of the Herman-Avila-Bochi formula for Lyapunov exponents of $\SL(2,\R)$-cocycles}

\author[A.~T.~Baraviera]{Alexandre T. Baraviera}
\address{Departamento de Matemática\\
Universidade Federal do Rio Grande do Sul\\
Porto Alegre, RS, Brasil}
\email{baravi@mat.ufrgs.br}

\author[J.~Lopes~Dias]{Jo\~ao Lopes Dias}
\address{Departamento de Matem\'atica and Cemapre, ISEG\\
Universidade T\'ecnica de Lisboa\\
Rua do Quelhas 6, 1200-781 Lisboa, Portugal}
\email{jldias@iseg.utl.pt}

\author[P.~Duarte]{Pedro Duarte}
\address{Departamento de Matemática and Cmaf \\
Faculdade de Ciências\\
Universidade de Lisboa\\
Campo Grande, Edifício C6, Piso 2\\
1749-016 Lisboa, Portugal 
}
\email{pduarte@ptmat.fc.ul.pt}

\date{\today}

\begin{abstract}
We study the geometry of the action of $\SL(2,\R)$ on the projective line in order to present a new and simpler proof of the Herman-Avila-Bochi formula.
This formula gives the average Lyapunov exponent of a class of 1-families of $\SL(2,\R)$-cocycles.

\end{abstract}

\maketitle


\section{Introduction}

A fundamental problem in smooth dynamics is the determination of the Lyapunov exponents of a given system.
These values correspond to the exponential rate of divergence or convergence of nearby orbits along prescribed directions. 
A positive Lyapunov exponent implies hyperbolic behaviour of orbits, which might produce very complicated dynamics. 
On the other hand, a negative Lyapunov exponent indicates that the orbits are fast converging and thus dynamics should be simpler.

Lyapunov exponents exist almost everywhere in phase space by the Oseledets theorem. However, their computation is typically a hard problem that has only been overcome by the use of numerical techniques. 
In fact, there are very few non-trivial examples outside uniform hyperbolicity for which their values (or even the signs) have been computed analytically.
Criteria for positive Lyapunov exponents for non-uniformly hyperbolic systems can be found in~\cite{K,M,W}.

In this paper we treat a remarkable example where the average Lyapunov exponents of families of cocycles can be explicitly computed. 
Herman~\cite{H} was the first to present it in the context of products of $\SL(2,\R)$-matrices over an ergodic transformation, and found a lower bound for the average upper Lyapunov exponent. 
Later, Avila and Bochi~\cite{AB} showed that Herman's lower bound was the actual value of the average exponent.
We present below the setting and results related to this problem, and give an alternative proof of the Herman-Avila-Bochi formula.
Our approach simplifies considerably the analysis by looking at simple geometric consequences of the action of the matrices on the projective line $\Pp^1$.

Let $(X,\mu)$ be a probability space, a measurable $\mu$-preserving ergodic transformation $f\colon X\to X$, and a $\mu$-integrable function $A\colon X\to\SL(2,\R)$.
We want to study the dynamics of the linear cocycle
$(f,A)\colon X\times\SL(2,\R)\to X\times\SL(2,\R)$ given by
$$
(f,A)(x,y)=(f(x),A(x)\,y).
$$
Its iterations are also linear cocycles
$$
(f,A)^n=(f^n,A_n),
$$
where 
$$
A_n(x)=A(f^{n-1}(x))\dots A(f(x)),
\qquad
n\in\N.
$$
Due to the above vector bundle structure we call the space $X$ the base, whilst $\SL(2,\R)$ is the fiber.

We deal with the question of obtaining the largest Lyapunov exponent on the fiber for the above cocyles. 
This is given by the asymptotic exponential growth of the norm of the product of matrices, measured by the {\it fiber upper Lyapunov exponent} of $(f,A)$, 
\begin{equation} 
\lambda(f,A)=\lim_{n\to+\infty} \frac{1}{n}\,\int_X \log\nrm{ A_n(x)  }\,d\mu(x).
\end{equation}

By considering the rotation by an angle $\theta$,
\begin{equation} \label{Rtheta}
R_\theta=\matriz{\cos \theta}{-\sin \theta}{\sin \theta}{\cos \theta},
\end{equation}
we focus on the 1-family of cocycles $\theta\mapsto(f,R_\theta A)$.
Using a sub-harmonicity ``trick'', Herman showed the following inequality for the average Lyapunov exponent inside this family.

\begin{teor}[Herman~\cite{H}]
\begin{equation*} 
\frac{1}{2\pi}\, \int_0^{2\pi} \lambda(f,R_\theta\,A)\,d\theta
\geq  \int_X \log\left(\frac{ \nrm{A(x)} + \nrm{A(x)}^{-1} }{2} \right) \,d\mu(x) \;.
\end{equation*}
\end{teor}

Roughly, Herman's method consists in showing that the function
$\theta\mapsto \lambda(f,R_\theta\,A)$ has a sub-harmonic extension to the
unit disk $\Dsk$ in the complex plane.
The inequality then follows from the sub-harmonicity property.
Later, under the same assumptions, Avila and Bochi improved Herman's inequality by showing that actually equality occurs.

\begin{teor}[Avila-Bochi~\cite{AB}]
\begin{equation*} 
\frac{1}{2\pi}\, \int_0^{2\pi} \lambda(f,R_\theta\,A)\,d\theta
= \int_X \log\left(\frac{ \nrm{A(x)} + \nrm{A(x)}^{-1} }{2} \right) \,d\mu(x) \;.
\end{equation*}
\end{teor}

As an example, this theorem applies immediately to the cocycle over an ergodic rotation $f$ on the circle $\R/\Z$ with 
$A(x)=R_x\left[\begin{smallmatrix}c&0\\0&c^{-1}
\end{smallmatrix}\right]$
and $c\not=1$.
We then have $R_\theta A(x)=A(x+\theta)$ and so $\lambda(f,R_\theta A)=\lambda(f,A)$ is constant in the family. 
The Herman-Avila-Bochi formula above gives $\lambda(f,A)=\log(c+c^{-1})-\log2>0$.
We remark that examples as this one are delicate since a $C^0$-generic $\SL(2,\R)$-cocycle is uniformly hyperbolic or it has zero Lyapunov exponent almost everywhere~\cite{B}.

Notice that for $A\in\SL(2,\R)$,
$$ \log\left(\frac{ \nrm{A} + \nrm{A}^{-1} }{2} \right) =
\int_{\Pp^1} \log \nrm{A\,p}\,dp $$
(see~\cite[Proposition 3]{AB}).
Using Birkhoff's ergodic theorem, 
Avila and Bochi reduce the proof of Theorem~2 to the following one (see ~\cite[Theorem 12]{AB}), where $\rho(A)$ stands for the logarithm of the spectral radius of $A$.

\begin{teor}[Avila-Bochi~\cite{AB}]\label{FHAB}
Given matrices $A_1,\ldots, A_n\in\SL(2,\R)$,
\begin{align*}
\frac{1}{2\pi}\, \int_{0}^{2\pi} \rho(R_\theta A_n \ldots  R_\theta A_1)\, d\theta 
&= \sum_{j=1}^n \int_{\Pp^1} \log\nrm{A_j\,p}\,dp \;.
\end{align*}
\end{teor}

To prove the above formula they show that the sub-harmonic extension of 
$$\theta\mapsto \rho(R_\theta A_n\ldots R_\theta A_1),
$$
as in Herman's trick, is in fact harmonic.
We present here an alternative proof of Theorem~\ref{FHAB} based on a simple change of variable argument, which exploits instead the geometry of the action of $\SL(2,\R)$ on the projective line $\Pp^1$.

In section~\ref{sec symm} we present some properties of the $\SL(2,\R)$-action on $\Pp^1$, and complete our proof of Theorem~\ref{FHAB} in section~\ref{sec mat seq}.

\section{Symmetries of Matrix Actions}\label{sec symm}

Consider the circle group $\Pp^1=\R/\pi\,\Z$ as a model of the real projective line and 
denote by $m$ the  normalized Haar measure on $\Pp^1$. 
Given $p\in\Pp^1$ let $\ell_p$ denote the line spanned by the vector
$$v_p=(\cos p,\sin p)\in\R^2.$$

For a matrix $A\in\SL(2,\R)$ the action $\Phi_A\colon\Pp^1\to\Pp^1$ of $A$ on $\Pp^1$ is characterized by the relation
\begin{equation}
\ell_{\Phi_{A}(p)}= A\ell_p,
\qquad
p\in\Pp^1.
\end{equation}
Its derivative is related to expansivity by
\begin{equation}
\Phi_A'(p)= \frac{1}{\nrm{A\,v_p}^2},
\qquad 
p\in\Pp^1.
\end{equation}
Moreover, define $H_A\colon\Pp^1\to \Pp^1$ as
$$
H_A(p)=p-\Phi_A(p).
$$

Notice that $\theta\mapsto R_\theta$ induces a well-defined map from
$\Pp^1$ to $\PSL(2,\R)$, where the rotation matrices $R_\theta$ are defined in~\eqref{Rtheta}.
The function $H_A$ can thus be characterized by the eigenspace relation
\begin{equation}\label{eigenspace}
R_{H_A(p)} A\,\ell_p = \ell_p,
\qquad 
p\in\Pp^1.
\end{equation}
Finally, take $\rho_A\colon\Pp^1\to \R$ to be a function which measures the expansivity of the action of $A$ as
$$
\rho_A(p)=\log \nrm{A\,v_p}.
$$
We then have the following properties.

\begin{prop}\label{PhiA:symmetries}
For every non-orthogonal matrix $A\in\SL(2,\R)$ there is
a unique analytic map $\Psi_A\colon\Pp^1\to\Pp^1$ such that:
\begin{enumerate}
\item $\Psi_A\circ\Psi_A = \mbox{id}_{\Pp^1}$,
\item $H_A\circ\Psi_A = H_A$,
\item $\rho_A\circ\Psi_A = -\rho_A$,
\item $\Psi_A'=-\Phi_A'$,
\item $H_A'=1+\Psi_A'$.
\end{enumerate}
\end{prop}

\begin{figure}[h]
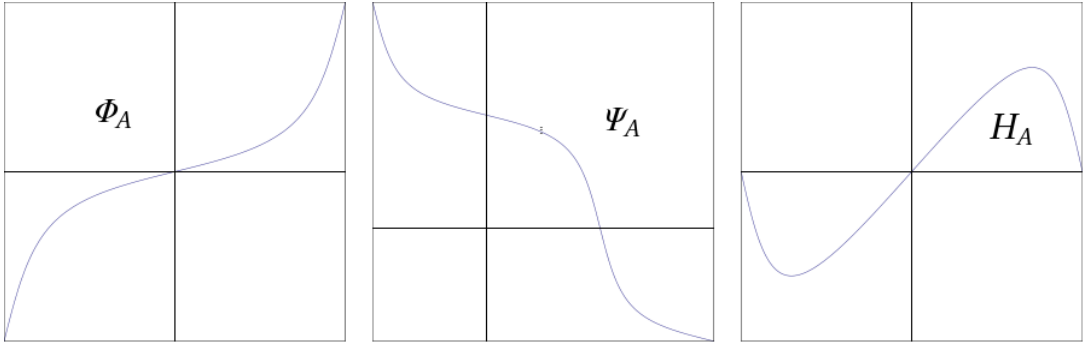

\figura{3g}{.45}
\caption{Functions  $\Phi_A$, $\Psi_A$  and $H_A$ }
\end{figure} 

\dem The uniqueness of such $\Psi_A$ is obvious since the pre-image $(H_A)^{-1}(p)$ of each
regular value $p\in H_A(\Pp^1)$ consists exactly of two points which must be inter-changed by $\Psi_A$.

By singular value decomposition, there exist $S,R\in\OO(2,\R)$ and $\lambda>1$ such that $A=S D R$, where
$$
D= \matriz{\lambda}{0}{0}{\lambda^{-1}}.
$$
Let $M = R^{-1}  K R$ with 
$$
K=\matriz{0}{\lambda^{-1}}{\lambda}{0}.
$$
We claim that $\Psi_A=\Phi_M$ is the required involution.

Since $K^2=I$, we have $M^2=I$ and item 1 follows. 
Notice also that
$$ \nrm{M\,v}=\nrm{K\,R\,v}=\nrm{D\,R\,v}=\nrm{A\,v}.$$
Hence,
$$ \nrm{A\,\Phi_M(v)} =
\frac{\nrm{A\,M\,v}}{\nrm{M\,v}} = \frac{\nrm{D\,K\,R\,v}}{\nrm{A\,v}} = 
\frac{\nrm{ R\,v}}{\nrm{A\,v}} = \frac{\nrm{ v}}{\nrm{A\,v}},$$
which proves item 3.

Next assume that $A$ is symmetric. We have $S=R^{-1}$, i.e.  $A=R^{-1}\,D \,R$,
and for this case
$$A\,M= R^{-1}\,D\,K\, R = R^{-1}\,\matriz{0}{1}{1}{0}\, R $$
is an isometric involution. 
Thus,
\begin{align*}
H_A\circ \Psi_A(p) & = \Phi_M(p) - \Phi_{A\,M}(p) = \Phi_{A\,M} (\Phi_{A\,M}(p)) - \Phi_{A\,M}(\Phi_M(p))\\
& = p - \Phi_{A\,M^2}(p) = p - \Phi_A(p) = H_A(p)\;.
\end{align*}
The general case, where  $A\neq R^{-1}\,D\,R$, now follows because $H_A-H_{R^{-1}\,D\,R}$ is a constant function. 
This implies that $H_A$ and $H_{R^{-1}\,D\,R}$ share the same involution
 $\Phi_M$.

As remarked above $\nrm{A\,v}=\nrm{M\,v}$.
Therefore,
$$ \Phi_A'(p) = \frac{1}{\nrm{A\,v_p}^2} = \frac{1}{\nrm{M\,v_p}^2} = -\Phi_M'(p)= -\Psi_A'(p).$$
Finally, item 5 follows from item 4 since $H_A' = 1-\Phi_A' = 1+\Psi_A'$.
\cqd

In our proof of the Herman-Avila-Bochi formula we will use the following abstract change of variables argument.

\begin{prop}\label{int:rho}
Consider an integrable function $\rho\colon I\to\R$ and a smooth involution
$\Psi\colon I\to I$ such that $\rho\circ \Psi=-\rho$. 
Then,
$$ \frac{1}{2}\,\int_I \rho(t)\,(1+\Psi'(t))\,dt =  \int_I \rho(t)\,dt \;.$$
\end{prop}

\dem
Let $I_{+}=\rho^{-1}(0,+\infty)$ and $I_{-}=\rho^{-1}(-\infty,0)$,  so that
$\Psi I_{+} = I_{-}$ and $\Psi I_{-} = I_{+}$. So, 
\begin{align*}
\int_I \rho(t)\,dt &= \int_{I_{+}} \rho(t)\,dt + \int_{I_{-}} \rho(t)\,dt\\
&= \int_{I_{+}} \rho(t)\,dt -  \int_{I_{+}} \rho\circ\Psi(t)\, \Psi'(t)\,dt\\
&= \int_{I_{+}} \rho(t)\,dt +  \int_{I_{+}} \rho(t)\, \Psi'(t)\,dt\\
&= \int_{I_{+}} \rho(t)\,(1+\Psi'(t) )\,dt \;. 
\end{align*}
Similarly, $\int_I \rho(t)\,dt = \int_{I_{-}} \rho(t)\,(1+\Psi'(t) )\,dt $.
Hence, the claim follows.
\cqd

Our next proposition is a special case of Theorem~\ref{FHAB}.
The proof illustrates how the previous argument applies.

\begin{prop}
For any matrix $A\in\SL(2,\R)$,
$$\frac{1}{2\pi}\, \int_0^{2\pi} \rho(R_\theta\,A)\,d\theta
= \int_{\Pp^1} \log \nrm{A\,p}\,dp.$$
\end{prop}

\dem
We will use the change of variable $\theta=H_A(p)$.
Notice first that $R_\theta\,A$ is elliptic iff $\theta$ lies outside the range of $H_A$, in which case the logarithm of the spectral radius of $R_\theta\, A$ is zero, i.e. $\rho(R_\theta\, A)=0$.
For the remaining values of $\theta$, the eigenspace property~(\ref{eigenspace}) implies that $\rho(R_{H_A(t)}\,A)=\Mod{\log \nrm{A\,v_t}	}$. 
Therefore,
\begin{align*}
 \frac{1}{\pi} \int_{-\pi/2}^{\pi/2}   \rho(R_\theta\,A)\,d\theta& =
\frac{1}{2\pi} \int_{-\pi/2}^{\pi/2} \Mod{\log \nrm{A\,v_t}}\, \Mod{H_A'(t)} \, dt \\ 
&=\frac{1}{2\pi} \int_{-\pi/2}^{\pi/2} \log \nrm{A\,v_t}\, H_A'(t) \, dt  \\
&=\frac{1}{2\pi} \int_{-\pi/2}^{\pi/2} \log \nrm{A\,v_t}\,\left(1+\Psi_A'(t) \right) \, dt  \\
&=   \int_{\Pp^1} \log \nrm{A\,p}\,dp \;.
\end{align*}
On the first step, the factor $1/2$ appears because the map $H_A$ covers twice the set of parameters
$\theta$ which correspond to  a real eigenvalue of $R_\theta\,A$.
The second equality follows because
$\log\nrm{A\,v_t}$ and $H_A'(t)= 1-\nrm{A\,v_t}^{-2}$ have the same sign for every $t$.
Then we use item 5 of Proposition~\ref{PhiA:symmetries},
and the final step is a consequence of Proposition~\ref{int:rho}.
\cqd

\section{Matrix Sequence Actions}\label{sec mat seq}

We call {\em matrix word} to any finite sequence of matrices
$$
\underline{A} = (A_1,\ldots, A_n)
$$
with $A_1,\ldots, A_n\in\SL(2,\R)$, and $n\in\N$ is the length of the word.
So, we denote by $\SL^n(2,\R)$ the space of all $\SL(2,\R)$-matrix words of length $n$.
For such a word we define the product
$$
R_\theta\underline{A} = (R_\theta\,A_1)\,(R_\theta\,A_2)\,\ldots\, (R_\theta\,A_n).
$$
Given any other matrix word $\underline{B} = (B_1,\ldots, B_k)$, we have 
$$
R_\theta(\underline{A}\,\underline{B}) = (R_\theta\underline{A})\,(R_\theta\underline{B}),
$$
where $\underline{A}\,\underline{B}$ stands for the concatenated word 
$(A_1,\ldots, A_n,B_1,\ldots, B_k)$.

Moreover, we choose the maps
$\Phi_{\underline{A}}:\Pp^1\times\Pp^1 \to \Pp^1$  
and $H_{\underline{A}}:\Pp^1\times\Pp^1 \to \Pp^1$
by 
$$
\Phi_{\underline{A}}(\theta,p) = \Phi_{R_\theta \underline{A}}(p)
\quad\text{and}\quad
H_{\underline{A}}(\theta,p) = p- \Phi_{\underline{A}}(\theta,p),
$$
respectively.

\begin{prop}
Given any  word $\underline{A} \in\SL^n(2,\R)$,
there exists $n$ analytic functions $\widetilde{H}_j(p) = \widetilde{H}_{\underline{A},j}(p)$, $j=1,\ldots, n$, implicitely 
defined by
$H_{\underline{A}}(\widetilde{H}_j(p), p ) = 0 $.
\end{prop}

\dem
The map  $\theta \mapsto H_{\underline{A}}(\theta, p )=p-\Phi_{\underline{A}}(\theta)$
is an expanding map of degree $n$. 
So, for each $p\in\Pp^1$
there are exactly $n$ points $\theta_j\in\Pp^1$, such that
$H_{\underline{A}}(\theta_j, p )=0$. By an implicit function theorem argument,
locally, each $\theta_j= \widetilde{H}_j(p)$ is an analytic function of $p$,
and we are left to prove that these local functions can be glued to form
$n$ global analytic functions. 
By defining $M$ as the union of these local manifolds, it is enough to prove that $M$ is a compact $1$-dimensional
manifold with $n$ connected components.

We can write $M$ as a pre-image 
$M=(G_{\underline{A}})^{-1}\left(\Pp^1\times\{0\}\right)$
of the map
$G_{\underline{A}}\colon \Pp^1\times \Pp^1\to\Pp^1\times\Pp^1$\, defined by \, $G_{\underline{A}}(\theta,p)=
\left( p, H_{\underline{A}}(\theta,p)\right)$. 
Its derivative is
$$
DG_{\underline{A}}(\theta,p) = \matriz{\frac{\partial H_{\underline{A}}}{\partial\theta} }{0}{\ast}{1},
$$
so $G_{\underline{A}}$ has no critical points and $M$ is a compact analytic $1$-dimensional manifold. 
Since $H_{\underline{A}}\colon\Pp^1\to\Pp^1$ is a map of zero degree, 
$G_{\underline{A}}$ induces a linear endomorphism on the homology space
$H_1(\Pp^1\times\Pp^1;\R)=\R^2$ whose action is given by the matrix
$$
\matriz{-n}{0}{0}{1}.
$$
Thus, $M=(G_{\underline{A}})^{-1}\left(\Pp^1\times\{0\}\right)$ must be the union of
$n$ homotopically non-trivial closed curves, which are precisely the graphs of
the functions $\widetilde{H}_{\underline{A},j}$.
\cqd

The functions $\widetilde{H}_{\underline{A},j}$ can also be characterized by the eigenspace relation 
\begin{equation}
R_{\widetilde{H}_{\underline{A},j}(p)}\underline{A}\,\ell_p = \ell_p 
\end{equation}
which implies that
\begin{equation*}
R_{\widetilde{H}_{\underline{A},j}(p)}\underline{A}\,v_p = \pm \nrm{R_{\widetilde{H}_{\underline{A},j}(p)}\underline{A}}\,v_p 
\end{equation*} 
That is, the matrix $R_{\theta}\underline{A}$ has the real eigenvector $v_p$ iff
$\theta=\widetilde{H}_{\underline{A},j}(p)$ for some $j=1,\ldots, n$.
This shows the following proposition.

\begin{prop}\label{elliptic}
The matrix $R_\theta\,\underline{A}$ is elliptic iff $\theta$ is not in the range of any of the functions
$\widetilde{H}_{\underline{A},j}$ with $j=1,\ldots, n$.
Moreover, $\rho(R_\theta\underline{A})= \Mod{\log\nrm{R_\theta\underline{A}\,v_p}}$
whenever $\theta=\widetilde{H}_{\underline{A},j}(p)$ for some $j=1,\ldots, n$
and $p\in\Pp^1$. Otherwise, $\rho(R_\theta\underline{A})=0$.
\end{prop}

\begin{figure}[h]
\figura{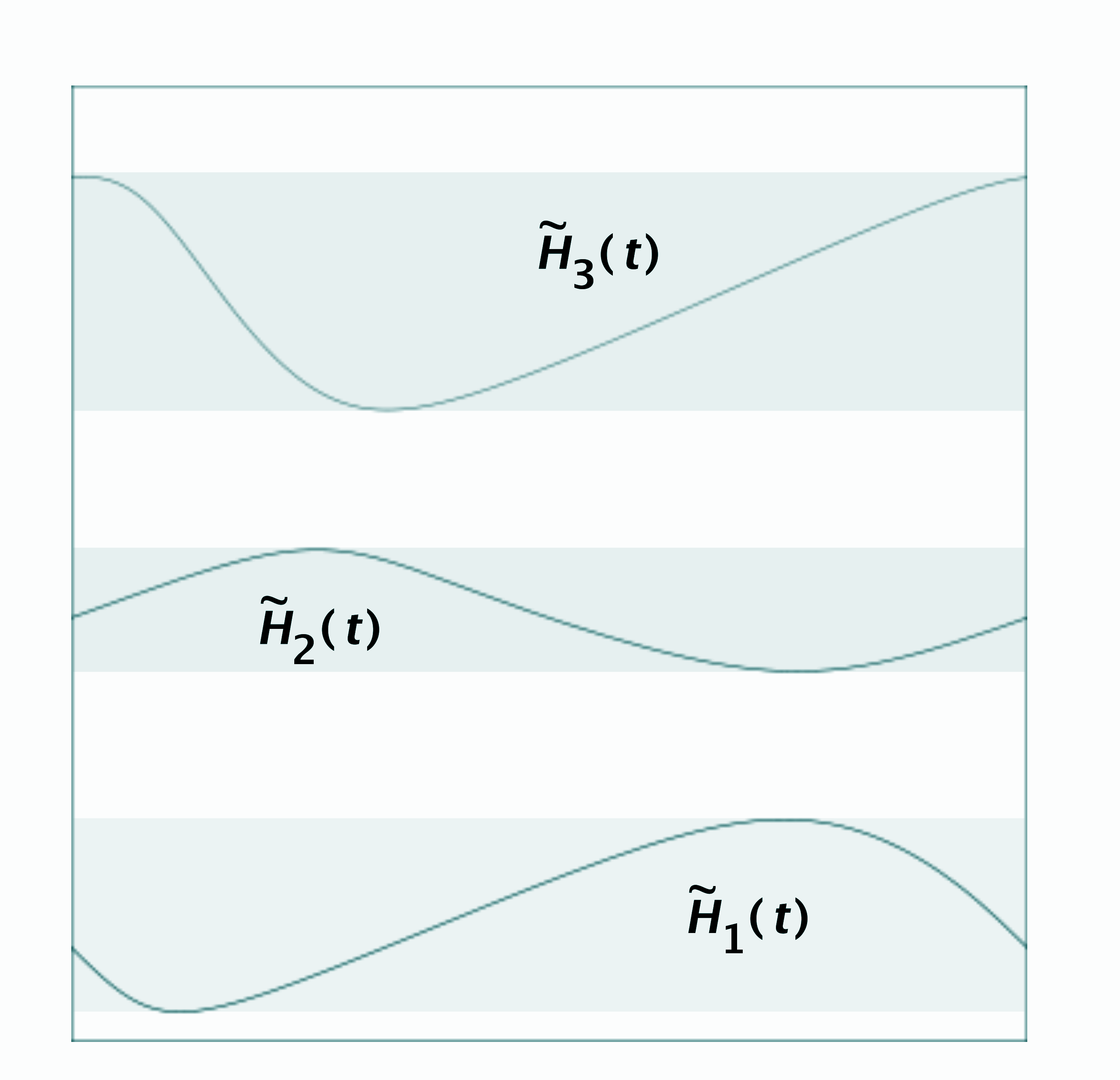}{.3}
\caption{Hyperbolic regions }
\end{figure}

Given $k,j=1,\ldots, n$ we define 
$$ v^k_{\underline{A},j}(p) := \frac{ R_{\widetilde{H}_j(p)}\,(A_{k+1},\ldots, A_n)\,v_p }
{ \nrm{ R_{\widetilde{H}_j(p)}\,(A_{k+1},\ldots, A_n)\,v_p }} \;, \quad \text{ and } $$
$$ \Phi^k_{\underline{A},j}(p) := \Phi_{(A_{k+1},\ldots, A_n)}(\widetilde{H}_j(p),p)  \;, $$\\
so that\, $v^k_{\underline{A},j}(p) = v_{\Phi^k_{\underline{A},j}(p)}$.
We also define
\begin{align*}
 J_{k}(\underline{A}) & := \frac{1}{2\pi}\, \sum_{j=1}^n 
\int_{-\pi/2}^{\pi/2} \log \nrm{ A_k\,v^k_{\underline{A},j}(p)   } \,\widetilde{H}_{\underline{A},j}'(p) \, dp \;.
\end{align*}

\begin{prop}\label{int:formula}
$$\frac{1}{2\pi}\, \int_{0}^{2\pi}  \rho(R_\theta\,\underline{A})\, d\theta  =
\sum_{k=1}^n J_k(\underline{A}). $$
\end{prop}

\dem
Let $I_j$ denote the range of $\widetilde{H}_{\underline{A},j}$.
Performing the change of variables $\theta=\widetilde{H}_{\underline{A},j}(p)$ in each interval
$I_j$, by Proposition~\ref{elliptic} we have
\begin{align*}
\frac{1}{2\pi}\, \int_{0}^{2\pi}  \rho(R_\theta\,\underline{A})\, d\theta  &=\frac{1}{\pi}\, \int_{-\pi/2}^{\pi/2}  \rho(R_\theta\,\underline{A})\, d\theta  =
 \frac{1}{\pi}\, \sum_{j=1}^n \int_{I_j}  \rho(R_\theta\,\underline{A})\, d\theta \\
&= \frac{1}{2\pi}\, \sum_{j=1}^n 
\int_{-\pi/2}^{\pi/2}  \rho(R_{\widetilde{H}_j(p)}\,\underline{A} )\, \Mod{ \widetilde{H}_j'(p)}\, dp \\
&= \frac{1}{2\pi}\, \sum_{j=1}^n 
\int_{-\pi/2}^{\pi/2} \Mod{\log \nrm{ R_{\widetilde{H}_j(p)}\,\underline{A} \,v_p}} \,\Mod{\widetilde{H}_j'(p)} \, dp \\
&= \frac{1}{2\pi}\, \sum_{j=1}^n 
\int_{-\pi/2}^{\pi/2} \log \nrm{ R_{\widetilde{H}_j(p)}\,\underline{A} \,v_p} \,\widetilde{H}_j'(p) \, dp \\
&= \frac{1}{2\pi}\, \sum_{j=1}^n \sum_{k=1}^n
\int_{-\pi/2}^{\pi/2} \log \nrm{ A_k\, \,v^k_{\underline{A},j}(p)} \,\widetilde{H}_j'(p) \, dp \\
&= \sum_{k=1}^n    J_k(\underline{A}).
\end{align*}
 The factor $1/2$ appears on the third step because the map $\widetilde{H}_{j}$ 
is a double cover of the interval  $I_j$.
The next step uses again Proposition~\ref{elliptic}.
Differentiating the relation
$\Phi_{\underline{A}}( \widetilde{H}_{j}(p), p ) = p$, which implicitely defines
$\widetilde{H}_{j}(p)$, we obtain
$$ \widetilde{H}_{j}'(p) = \frac{ 1- \frac{\partial \Phi_{\underline{A}} }{\partial p}( \theta,p) }{ \frac{\partial \Phi_{\underline{A}} }{\partial \theta}(\theta,p) }
= \frac{ 1- \frac{1}{\nrm{R_\theta\underline{A}\,v_p}^2} }{ \frac{\partial \Phi_{\underline{A}} }{\partial \theta}(\theta,p) } 
$$
with $\theta= \widetilde{H}_{j}(p)$.
Hence, since $\frac{\partial \Phi_{\underline{A}} }{\partial \theta}>0$, the numbers
$\widetilde{H}_{j}'(p)$, 
$$
1- {\nrm{R_{\widetilde{H}_{j}(p)} \underline{A}\,v_p}^{-2}}
\quad\text{and}\quad
\log \nrm{R_{\widetilde{H}_{j}(p)} \underline{A}\,v_p}
$$ 
have the same sign, which explains the fifth step. Step six follows by cocycle additivity,
and by exchanging the two summations, we complete the proof.
\cqd

\begin{lema}\label{measpres}
Given a matrix word  $\underline{A} \in\SL^n(2,\R)$ and $p\in\Pp^1$, the map 
$f\colon\Pp^1\to\Pp^1$\,
$f(\theta)= \Phi_{\underline{A}}(\theta ,p)$ \, 
is  an expanding map on $\Pp^1$ with degree $n$ which preserves the Haar measure on $\Pp^1$.
\end{lema}

\dem Denote by $m$ the Haar measure, both on $\Pp^1$, and on
$\Su^1=\{\, z\in\C\,:\,\Mod{z}=1\,\}$.
Let $\Dsk=\{\,z\in\C\colon \Mod{z}<1\,\}$ be the unit disk.
For each matrix 
$$A=\matriz{a}{b}{c}{d}\in\SL(2,\R),$$ 
define $M_A\colon\R\cup\{\infty\}\to\R\cup\{\infty\}$ by
$M_A(x)= \frac{ a\,x+b}{c\,x + d}$. Then $\xi\circ \Phi_A = M_A\circ \xi$, where $\xi\colon\Pp^1\to \R\cup\{\infty\}$
is the map $\xi(p) = \tan p $. 

Consider the M\"{o}bius transformation 
$\eta(z)=\frac{1+i\,z}{1-i\,z}$, which maps $\R\cup\{\infty\}$ onto $\Su^1$ and
let $\psi= \eta\circ \xi$. 
The fundamental formula of trigonometry implies
that $\psi(p)=e^{2 i p}$.
Notice that $\psi$ is a continuous group isomorphism. Hence $\psi_\ast m = m$.

Define now $\widehat{M}_A\colon \Su^1\to \Su^1$,
$$
\widehat{M}_A(z) = \psi\circ \Phi_A\circ \psi^{-1}(z)= \eta\circ M_A\circ \eta^{-1}(z),
$$
which extends to a M\"{o}bius transformation on the Riemann sphere that preserves the circle $\Su^1$.
The linear fractional map $\widehat{M}_A(z)$ satisfies the symmetry relation
 
$$
\widehat{M}_A(\obar{z}^{-1}) = \obar{\widehat{M}_A(z)}^{-1}.
$$ 
It has a single zero inside the disk $\Dsk$, and a single pole outside.
With this notation,  define $\widehat{f}\colon \Su^1\to \Su^1$ 
$$ 
\widehat{f}(z) = z\, \widehat{M}_{A_1}\left( z\,\widehat{M}_{A_2}( \, \ldots\, z\,\widehat{M}_{A_n}(\psi(p)) \ldots ) \right).
$$
So, $\widehat{f}=\psi\circ f \circ \psi^{-1}$. Now, $\widehat{f}(z)$  is a rational function satisfying the symmetry relation
$$
\widehat{f}(\obar{z}^{-1}) = \obar{\widehat{f}(z)}^{-1},
$$
with  zeros inside the disk $\Dsk$, and poles outside\footnote{
Functions with these properties are finite Blaschke products, see~\cite{Con}.}.
The map $\widehat{f}$ is analytic on $\Dsk$ with $\widehat{f}(0)=0$. We claim that
this property implies that $\widehat{f}_\ast m = m$,
which in turn will imply  $f_\ast m = m$, and finish the proof.
To see this, take any continuous function $\varphi\colon \Su^1\to\C$. By the Dirichelet principle
this function has a continuous extension $\widetilde{\varphi}\colon \obar{\Dsk}\to\C$
\, which is harmonic on $\Dsk$. We refer it as the {\em harmonic} extension of $\varphi$.
Since $\widehat{f}$ is analytic on $\Dsk$,\, $\widetilde{\varphi}\circ\widehat{f}$
is the harmonic extension of $\varphi\circ \widehat{f}$. Therefore, by the Poisson formula\,
$$ \int_{\Su^1} \varphi\circ \widehat{f}\, dm = \widetilde{\varphi}( \widehat{f} (0) ) =
\widetilde{\varphi}( 0 ) = \int_{\Su^1} \varphi \, dm\;,$$
which implies that $\widehat{f}_\ast m = m$.
\cqd

\begin{prop} \label{Hpme} For $\underline{A}=(A_1,\ldots, A_n)\in\SL^n(2,\R)$ and $p\in\Pp^1$,
$$\sum_{j=1}^n \widetilde{H}_{\underline{A},j}'(p)  = 1-\Phi_{A_n}'(p)  
= 1+ \Psi_{A_n}'(p) \;.$$
\end{prop}

\dem
We have for the matrix word  $\underline{B}=(A_{n-1}^{-1},A_{n-2}^{-1},\ldots, A_1^{-1},I)$,
$$\Phi_{\underline{B}}(\theta,p) = \theta + \Phi_{A_{n-1}}^{-1}\left( \theta 
+ \Phi_{A_{n-2}}^{-1}( \,\ldots \,\Phi_{A_1}^{-1}(\theta + p )\,\ldots\, )\,\right) $$
and hence
\begin{equation}\label{AB}
 \Phi_{\underline{B}}\left( -\theta, \Phi_{\underline{A}}(\theta, p ) \right) = \Phi_{A_n}(p)  
\end{equation}
Differentiating~(\ref{AB})  w.r.t. $\theta$ and $p$ we get respectively 
$$ -\frac{\partial \Phi_{\underline{B}}}{\partial \theta}\left( -\theta, \Phi_{\underline{A}}(\theta, p ) \right) +
\frac{\partial \Phi_{\underline{B}}}{\partial p}\left( -\theta, \Phi_{\underline{A}}(\theta, p ) \right) \,
\frac{\partial \Phi_{\underline{A}}}{\partial \theta}\left(\theta, p \right) =0 \;, $$
$$  \frac{\partial \Phi_{\underline{B}}}{\partial p }\left( -\theta, \Phi_{\underline{A}}(\theta, p ) \right) \,
\frac{\partial \Phi_{\underline{A}}}{\partial p}\left(\theta, p \right) =\Phi_{A_n}'(p) \;. $$
Hence
\begin{equation}\label{DPhi}
 \frac{ \frac{\partial \Phi_{\underline{A}}}{\partial p}\left(\theta, p \right) }{
\frac{\partial \Phi_{\underline{A}}}{\partial \theta }\left(\theta, p \right) } =
\frac{ \Phi_{A_n}'(p)  }{
\frac{\partial \Phi_{\underline{B}}}{\partial \theta }\left(-\theta, \Phi_{\underline{A}}(\theta, p ) \right) } 
\end{equation}
Write $\widetilde{H}_j = \widetilde{H}_{\underline{A},j}$. Differentiating the defining relation
$\Phi_{\underline{A}}( \widetilde{H}_j(p), p ) = p$ and writing
$\theta_j=\widetilde{H}_j(p) $ for $j=1,\ldots, n$  we obtain by~(\ref{DPhi}) that
\begin{align*}
\sum_{j=1}^n \widetilde{H}_j'(p) &=
\sum_{j=1}^n \frac{ 1- \frac{\partial \Phi_{\underline{A}}}{\partial p}(\theta_j, p)}{
\frac{\partial \Phi_{\underline{A}}}{\partial \theta}(\theta_j, p) }\\
&=
\sum_{j=1}^n \frac{1}{
\frac{\partial \Phi_{\underline{A}}}{\partial \theta}(\theta_j, p) } -
\sum_{j=1}^n \frac{\Phi_{A_n}'(p)}{
\frac{\partial \Phi_{\underline{B}}}{\partial \theta}(-\theta_j, \Phi_{\underline{A}}( \theta_j, p ) ) }
 = 1-\Phi_{A_n}'(p).
\end{align*}

By Lemma~\ref{measpres},  
since the $n$ points $\theta_j$ are the pre-images of $p$ by the  measure preserving expanding map
$\theta\mapsto \Phi_{\underline{A}}(\theta,p)$, 
$$\sum_{j=1}^n \frac{1}{
\frac{\partial \Phi_{\underline{A}}}{\partial \theta}(\theta_j, p) } = 1 \;. $$ 
Similarly,
$$\sum_{j=1}^n \frac{1}{
\frac{\partial \Phi_{\underline{B}}}{\partial \theta}(-\theta_j, \Phi_{\underline{A}}( \theta_j, p ) ) } =
\sum_{j=1}^n \frac{1}{
\frac{\partial \Phi_{\underline{B}}}{\partial \theta}(-\theta_j, p ) } = 1$$
because the $n$ points $-\theta_j$ are the pre-images of $\Phi_{A_n}(p)$ by the  measure preserving expanding map
$\theta\mapsto \Phi_{\underline{B}}(\theta,p)$.
\cqd

\begin{prop} \label{cyclic}  For each $k=1,\ldots, n$,
$$ J_k(A_1,\ldots, A_n) = J_n(A_{k+1},\ldots, A_n, A_1, \ldots, A_k) $$
\end{prop}

\dem By definition we have $\Phi_{\underline{A}}(\widetilde{H}_{\underline{A},j}(p),p) = p$.
Setting 
$$\underline{B}=(A_{k+1},\ldots, A_n, A_1, \ldots, A_k)\; $$
since the matrices
$R_{\widetilde{H}_j(p)}\underline{A}$ and $R_{\widetilde{H}_j(p)}\underline{B}$ are conjugate
by $R_{\widetilde{H}_j(p)}(A_{k+1},\ldots, A_n)$ we get
$$ \Phi_{\underline{B}}(\widetilde{H}_{\underline{A},j}(p), \Phi^k_{\underline{A},j}(p) ) =  \Phi^k_{\underline{A},j}(p)\;,$$
and hence, for $j=1,\ldots, n$,
$$ \widetilde{H}_{\underline{A},j}(p) = \widetilde{H}_{\underline{B},j}\left( \Phi^k_{\underline{A},j}(p)\right)\;.$$
Differentiating this relation we obtain 
\begin{align*}
J_k(\underline{A}) &= \frac{1}{2\pi}\,\sum_{j=1}^n \int_{-\pi/2}^{\pi/2}
\log \nrm{ A_k\,v^k_{\underline{A},j}(p)   } \,\widetilde{H}_{\underline{A},j}'(p) \, dp \\
&= \frac{1}{2\pi}\,\sum_{j=1}^n \int_{-\pi/2}^{\pi/2}
\log \nrm{ A_k\,v_{\Phi^k_{\underline{A},j}(p)}   } \,
\widetilde{H}_{\underline{B},j}'\left(\Phi^k_{\underline{A},j}(p)\right)\,  \,
(\Phi^k_{\underline{A},j})'(p)\, dp \\
&= \frac{1}{2\pi}\,\sum_{j=1}^n \int_{-\pi/2}^{\pi/2}
\log \nrm{ A_k\,v_{p}   } \,
\widetilde{H}_{\underline{B},j}'(p)\, dp  = J_n(\underline{B})\;.
\end{align*}
\cqd

\begin{prop} \label{Jk} For each $k=1,\ldots, n$,
$$ J_k(\underline{A}) = \int_{\Pp^1} \log \nrm{A_k\,p}\, dp $$
\end{prop}

\dem
By Proposition~\ref{cyclic} it is enough to consider the case $k=n$.
Then combining propositions~\ref{Hpme} and~\ref{int:rho}, we get the third
and fourth equalities below
\begin{align*}
J_n(\underline{A}) &= \frac{1}{2\pi}\,\sum_{j=1}^n \int_{-\pi/2}^{\pi/2}
\log \nrm{ A_n\,v_{p}   } \,
\widetilde{H}_{\underline{A},j}'(p)\, dp \\
&= \frac{1}{2\pi}\, \int_{-\pi/2}^{\pi/2}
\log \nrm{ A_n\,v_{p}   } \,\left( \sum_{j=1}^n
\widetilde{H}_{\underline{A},j}'(p) \right)\, dp \\
&= \frac{1}{2\pi}\, \int_{-\pi/2}^{\pi/2}
\log \nrm{ A_n\,v_{p}   } \,\left( 1+\Psi_{A_n}'(p)  \right)\, dp \\
&= \frac{1}{\pi}\, \int_{-\pi/2}^{\pi/2}
\log \nrm{ A_n\,v_{p}   } \,  dp  = \int_{\Pp^1} \log \nrm{ A_n\,v_{p}   } \,  dp\;.
\end{align*}
\cqd

Theorem~3 is a corollary of Propositions~\ref{int:formula} and~\ref{Jk}.

\section*{Acknowledgements}

This work was partially supported by Funda\c c\~ao para a Ci\^encia e a Tecnologia through the project ``Randomness in Deterministic Dynamical Systems and Applications'' ref. PTDC/MAT/105448/2008.


\thispagestyle{empty}

\end{document}